\documentclass[a4paper,12pt]{article}

\usepackage{amsmath,amsfonts,amssymb,graphicx,color,array}

\newtheorem{newremark}{Remark}[section]

\newtheorem{theorem}{Theorem}

\newtheorem{example}[theorem]{Example}
\newtheorem{lemma}[theorem]{Lemma}

\newcommand{\I}{\mathbb{I}}

\newcommand{\e}{{\rm e}}
\newcommand{\R}{\mathbb{R}}

\newcommand{\N}{\mathbb{N}}

\newcommand{\cF}{ {\cal F} }
\newcommand{\cE}{ {\cal E} }
\newcommand{\cA}{ {\cal A} }

\newcommand{\cor}[1]{{#1}}

\begin{document}

\title{Fast convergent method for the $m$-point problem in Banach space}

\author{D.O.~Sytnyk\thanks{Institute of
Mathematics of NAS of Ukraine, 3 Tereshchenkivs'ka Str., Kyiv-4, 01601, Ukraine
({\tt sytnik@imath.kiev.ua}).} \and V.B.~Vasylyk\thanks{Institute of
Mathematics of NAS of Ukraine,
 3 Tereshchenkivs'ka Str., Kyiv-4, 01601, Ukraine
({\tt vasylyk@imath.kiev.ua}).}
 }

\date{\null}
\maketitle

\begin{abstract}
The $m$-point nonlocal  problem for the first order differential
equation with an operator coefficient in a Banach space $X$ is
considered. An exponentially convergent algorithm is proposed and
justified provided that the operator coefficient $A$
 is strongly positive and some existence and uniqueness conditions are fulfilled.
 This algorithm  is based on representations of operator functions by a Dunford-Cauchy integral along a hyperbola
 enveloping the spectrum of $A$ and on the proper quadratures involving short sums of resolvents. The efficiency of the
 proposed algorithms is demonstrated by numerical examples.
\end{abstract}



\section{Introduction}

The m-point initial (nonlocal) problem for a differential equation with the nonlocal  condition $u(t_0) + g(t_1;\dots ;
t_p; u) = u_0$ and a given function $g$ on a given point set  $P=\{0= t_0 < t_1 <\dots < t_p\}$ , is one of the
important topics in the study of differential equations. Interest in such problems {originates} mainly from some
physical problems with a control of the solution at $P$. For example, when the function $g(t_1;\dots ; t_p; u)$ is
linear, we will have the periodic problem $u(t_0)=u(t_1)$. Problems with nonlocal conditions arise in the theory of
physics of plasma \cite{samde}, nuclear physics \cite{LeunCh}, mathematical chemistry \cite{LeuOrt}, waveguides
\cite{Gord} etc.

Differential equations with operator coefficients in some Hilbert or Banach space can be considered as meta-models for
systems of partial or ordinary differential equations and are suitable for investigations using tools of the functional
analysis (see e.g. \cite{clem, krein}). Nonlocal problems can also be considered within this framework
\cite{Byszewski2, Byszewski3}.

In this paper we consider the following nonlocal $m$-point problem:
\begin{equation}\label{pr-st}
  \begin{split}
    &u^{\prime}_t+Au=f(t), \quad t \in [0,T]\\
    &u(0)+\sum_{k=1}^{m}\alpha_k u(t_k) =u_0,\quad 0<t_1<t_2<\ldots < t_m\leq T,
  \end{split}
\end{equation}
where $\alpha_k \in \R,$ $k=\overline{1,m},$ $f(t)$ is a given
vector-valued function with values in \cor{a} Banach space $X,$ $u_0\in X.$
The operator $A$ with the domain $D(A)$ in a Banach space $X$ is
assumed to be a densely defined strongly positive (sectorial)
operator, i.e. its spectrum $\Sigma(A)$ lies in a sector of the
right half-plane with the vertex at the origin and and the resolvent
decays inversely proportional to $|z|$ at the infinity (see estimate
(\ref{estrez}) below).

Discretization methods for differential equations in Banach and
Hilbert spaces were intensively studied in the last decade (see e.g.
\cite{fernandez-lubich-palencia-schaedle, ghk, ghk3, gm2,
lopez-fernandez2, lopez-fernandez1, SSloanTho:98, sst1, vas1, vas2}
and the references therein). Methods from \cite{ ghk, ghk3,  gm2,
lopez-fernandez2, lopez-fernandez1, sst1, vas1, vas2} possess an
\cor{exponential convergence rate}, i.e. the error estimate in {an
appropriate} norm is of the type ${\cal O}(\e^{-N^\alpha}),$
$\alpha>0$ with respect to a discretization parameter $N \to
\infty$. For a given tolerance $\varepsilon$ such discretizations
provide optimal or nearly optimal computational complexity
\cite{ghk}. One of the possible ways to obtain exponentially
convergent approximations to {abstract differential equations} is
based on a representation of the solution through the Dunford
-Cauchy integral along a parametrized path enveloping the spectrum
of the operator coefficient. {Choosing a proper quadrature for this
integral we obtain a short sum of resolvents. Since the treatment of
such resolvents is usually \cor{the} most time consuming part of any
approximation this leads to a low-cost naturally parallelizeable
algorithms. } Parameters of the algorithms from \cite{ gm2,
lopez-fernandez1, sst1} were optimized in \cite{weid1, weid2} to
improve the convergence rate.

The aim of this paper is to construct an exponentially convergent
approximation to the problem for a differential equation with
$m$-point nonlocal condition in abstract setting (\ref{pr-st}). The
paper is organized as follows. { In Section \ref{stet} we discuss
the existence and uniqueness of the solution as well as its
representation through input data.}  An algorithm for the
homogeneous problem (\ref{pr-st}) is proposed in section
\ref{s-hom-pr}. The main result of this section \cor{is theorem}
\ref{hom-conv} about the exponential convergence rate of the
proposed discretization. The next section \ref{s-inh-pr} is devoted
to the exponentially convergent discretization of the inhomogeneous
problem. In section \ref{num-ex} we represent some numerical
examples which {confirm theoretical results from the previous
sections.}



\section{Existence and representation of the solution}\label{stet}

Let the operator $A$ in (\ref{pr-st}) be a densely defined strongly
positive (sectorial) operator in a Banach space $X$ with the domain
$D(A),$ i.e. its spectrum $\Sigma(A)$ lies in the sector (spectral
angle)
\begin{equation}\label{spectrA}
\Sigma = \left\{ z=\rho_0+r \e^{i\theta}:\quad r \in [0,\infty), \ \left|\theta\right|< \varphi < \frac{\pi}{2}
\right\}.
\end{equation}
Additionally outside the sector and on its boundary $\Gamma_\Sigma$ the following estimate for the resolvent holds true
\begin{equation}\label{estrez}
\left\|(zI-A)^{-1}\right\|\leq \frac{M}{1+\left|z\right|}.
\end{equation}
{The numbers $\rho_0,$ $\varphi$ are called } the spectral
characteristics of $A$.

The hyperbola
\begin{equation}\label{sp-hyp}
\Gamma_0=\{z(\xi)=\rho_0 \cosh{\xi}-i b_0 \sinh{\xi}: \; \xi \in
(-\infty, \infty ), \; b_0=\rho_0 \tan{\varphi}\}
\end{equation}
{in turn is \cor{referred} as a spectral hyperbola. It has a vertex at
$(\rho_0,0)$ and asymptotes which are parallel to the rays of the
spectral angle $\Sigma$.}

{A convenient representation of operator functions is the one
through the Dunford-Cauchy integral (see e.g. \cite{clem, krein})
\cor{where the integration path plays an important role}.  We choose the
following hyperbola }
\begin{equation}\label{int-hyp}
\Gamma_I=\{z(\xi)=a_I \cosh{\xi}-i b_I \sinh{\xi}: \; \xi \in
(-\infty, \infty )\},
\end{equation}
as the integration contour which envelopes the spectrum of $A$.

For an arbitrary vector $u_0 \in D(A^{m+1})$ the next equality holds
\begin{equation}\label{fp3}
\begin{split}
&\sum_{k=1}^{m+1}\frac{A^{k-1}u_0}{z^k}+\frac{1}{z^{m+1}}(zI-A)^{-1}A^{m+1}u_0=(zI-A)^{-1}u_0.
\end{split}
\end{equation}
This formula, together with
\begin{equation}\label{fp4}
  A^{-(m+1)}v=\frac{1}{2 \pi i}\int_{\Gamma_I}z^{-(m+1)}(zI-A)^{-1}v dz ,
\end{equation}
by setting $v=A^{m+1}u_0,$ yields the following representation
\begin{equation}\label{fp5}
\begin{split}
 & u_0=A^{-(m+1)}A^{m+1}u_0=\frac{1}{2 \pi
  i}\int_{\Gamma_I}z^{-(m+1)}(zI-A)^{-1}A^{m+1}u_0dz\\
  &=\int_{\Gamma_I}\left[(zI-A)^{-1}-\sum_{k=1}^{m+1}\frac{A^{k-1}}{z^k}\right]u_0
  dz.
\end{split}
\end{equation}

Various useful  properties of this representation as well as the next result were discussed in \cite{gm5}.
\begin{theorem}\label{t1}
Let $u_0 \in D(A^{m+\alpha})$ for some $m \in \N$ and $\alpha \in
[0,1],$ then the following estimate holds true:
\begin{equation}\label{fp11}
\begin{split}
 \left\|\left[(zI-A)^{-1}-\sum_{k=1}^{m+1}\frac{A^{k-1}}{z^k} \right]u_0\right\| &\le
\frac{1}{|z|^{m+1}} \frac{(1+M)K}{(1+|z|)^{\alpha}}\|A^{m+\alpha} u_0\|,\\
  &\forall \alpha \in[0,1], \; u_0 \in D(A^{m+\alpha}),
 \end{split}
\end{equation}
where the constant $K$ depends on $\alpha$ and $M$ only.
\end{theorem}

In \cite{Byszewski2, Byszewski3} it was proven that the solution of
the problem (\ref{pr-st}) exists and is unique provided that one of
the following two conditions is fulfilled:
\begin{equation}\label{um1}
\sum_{i=1}^{m} |\alpha_i| < 1,
\end{equation}
or
\begin{equation}\label{um2}
\sum_{i=1}^{m} |\alpha_i|e^{-\rho_0 t_i} < 1,
\end{equation}
with $\rho_0$  from \cor{(\ref{spectrA})}.

According to the Hille-Yosida-Phillips theorem \cite{yosida} the strongly positive
operator $A$ generates a one parameter semigroup $T(t)= e^{-t A}$
and the solution of (\ref{pr-st}) can be represented by
\begin{equation}\label{bp1int}
  u(t)=\e^{-At}u(0)+\int_0^t{\e^{-A(t-\tau)}f(\tau)}d\tau .
\end{equation}

Combining the nonlocal condition from (\ref{pr-st}) and (\ref{bp1int}) we obtain
\begin{equation}\label{bp1LS}
  \begin{split}
    u(t_i)=&\e^{-At_i}\left[u_0-\sum_{k=1}^{n}\alpha_k u(t_k) \right]+\int_0^{t_i}{\e^{-A(t_i-\tau)}f(\tau)}d\tau ,\\
    &i=\overline{1,m}.
  \end{split}
\end{equation}
Multiplying these equations by $\alpha_i$ and summing up over
$i=\overline{1,m}$ one gets

\begin{equation}\label{sum_t_i}
\begin{split}
  \sum_{i=1}^m \alpha_i u(t_i) = & \sum_{i=1}^m \alpha_i \e^{-At_i} u_0-\sum_{i=1}^m \alpha_i \e^{-At_i} \sum_{k=1}^m \alpha_k u(t_k) \\
  &+ \sum_{i=1}^m\alpha_i \int_0^{t_i} \e^{-A (t_i-\tau)}f(\tau) d\tau.
 \end{split}
\end{equation}

 Let us denote $$W=\sum_{i=1}^m \alpha_i u(t_i),$$ then
 (\ref{sum_t_i}) implies
 $$
 W = -\sum_{i=1}^m \alpha_i \e^{-At_i} W +\sum_{i=1}^m \alpha_i \e^{-At_i} u_0+ \sum_{i=1}^m \alpha_i \int_0^{t_i} \e^{-A (t_i-\tau)}f(\tau) d\tau,
 $$
Setting \cor{$B=\I +\sum\limits_{i=1}^m \alpha_i \e^{-At_i},$} we obtain
 $$
 B W = B u_0 -u_0 +\sum_{i=1}^m \alpha_i \int_0^{t_i} \e^{-A (t_i-\tau)}f(\tau) d\tau,
 $$
where $\I$ is an identity operator. Condition (\ref{um1}) or (\ref{um2})  guarantee the existence and boundedness of
$B^{-1}(A)$, therefore we  have
$$
 W = u_0- B^{-1} u_0 + B^{-1}\sum_{i=1}^m \alpha_i \int_0^{t_i} \e^{-A (t_i-\tau)}f(\tau) d\tau.
$$
Now,  the equations (\ref{bp1LS}) yield the  representation for the solution of (\ref{pr-st}) in the form
$u(t)=u_h(t)+u_{ih}(t)$, where $u_h(t)=\e^{-At} B^{-1} u_0$ is the solution of the homogeneous problem with the initial
condition $u_0$ and
 \begin{equation}\label{bp1IntRed}
 \begin{split}
   &u_{ih}(t)=-\e^{-At}  B^{-1}\sum_{i=1}^m \alpha_i \int_0^{t_i} \e^{-A (t_i-\tau)} f(\tau) d\tau +\int_0^t{e^{-A(t-\tau)}f(\tau)}d\tau
 \end{split}
 \end{equation}
is the solution of the inhomogeneous problem with zero initial condition.


\section{Homogeneous nonlocal problem}\label{s-hom-pr}

First of all we consider the solution $u_h(t)=\e^{-At} B^{-1}(A) u_0$ \cor{of the homogeneous} problem (\ref{pr-st}). Our aim in
this section is to construct an exponentially convergent method for its approximation.

Using the Dunford-Cauchy representation of $u_h(t)$ analogously to
\cite{gm5} we obtain
\begin{equation}\label{rep-hom0}
\begin{split}
 u_h(t)&=\frac{1}{2 \pi i}\int_{\Gamma_I} \e^{-zt}B^{-1}(z)(zI-A)^{-1} u_0 dz \\
 &=\frac{1}{2 \pi i}\int_{\Gamma_I} \frac{\e^{-zt}}{1 +\sum\limits_{i=1}^n \alpha_i \e^{-zt_i}} (zI-A)^{-1} u_0 dz.
\end{split}
\end{equation}
{Representation  (\ref{rep-hom0}) makes sense only when the function $\e^{-zt}B^{-1}(z)$ is analytic in the region
enveloped by $\Gamma_I.$} Let us show, that conditions like (\ref{um1}) or (\ref{um2}) guaranty this analyticity
\cite{krein}.

Actually, the analyticity of $\e^{-zt}B^{-1}(z)$ might only be vi\-olated when \\ $\sum\limits_{i=1}^n \alpha_i
\e^{-zt_i}=-1$, since in this case $B^{-1}(z)$ becomes unbounded. It is easy to see that for an arbitrary $z$ we have

\[
  \left|1+\underset{k=1}{\overset{m}{\sum}}\alpha_k e^{-zt_k}\right|
\ge \left|1-\underset{k=1}{\overset{m}{\sum}}|\alpha_k| \cdot
|e^{t_k  a_I \cosh{\xi} }|\cdot |e^{-i
 t_k b_I \sinh{\xi}}|\right|\]
\[\ge
\left|1-\underset{k=1}{\overset{m}{\sum}}|\alpha_k|  e^{-\rho_0 t_k}\right|>0 ,
\]
provided that (\ref{um2}) holds true.

 Using (\ref{fp5}) with $m=0$ we can modify the representation of  $u_h(t)$ as follows:
\begin{equation}\label{rep-hom}
\begin{split}
 u_h(t)&=\frac{1}{2 \pi i}\int_{\Gamma_I} \e^{-zt}B^{-1}(z)\left[(zI-A)^{-1}-\frac{1}{z}I \right]u_0 dz \\
 &=\frac{1}{2 \pi i}\int_{\Gamma_I} \frac{\e^{-zt}}{1 +\sum\limits_{i=1}^n \alpha_i \e^{-zt_i}} \left[(zI-A)^{-1}-\frac{1}{z}I \right]u_0
 dz.
\end{split}
\end{equation}
After discretization of the integral such  modified resolvent provides better convergence \cor{speed than} (\ref{rep-hom0})
in a \cor{neighbourhood} of $t=0$ (see \cite{gm5} for details).

Parameterizing the integral (\ref{rep-hom}) by (\ref{int-hyp}) we get
\begin{equation}\label{par-int}
 u_h(t)= \frac{1}{2 \pi i}\int_{-\infty}^{\infty}\cF(t,\xi) d \xi,
\end{equation}
with
\[
\cF(t,\xi)=F_A(t,\xi)u_0,
\]
\[
 F_A(t,\xi)=\frac{\e^{-z(\xi)t}z'(\xi)}{1 +\sum\limits_{i=1}^n \alpha_i \e^{-z(\xi)t_i}} \left[(z(\xi)I-A)^{-1}-\frac{1}{z(\xi)}I \right],
\]
\[
z'(\xi)= a_I \sinh{\xi}-i b_I \cosh{\xi}.
\]

Supposing $u_0 \in D(A^{\alpha}), \; 0<\alpha<1$ it was shown in \cite{gm5} that
\[
\begin{split}
  &\| \e^{-z(\xi)t}z'(\xi) \left[(z(\xi)I-A)^{-1}-\frac{1}{z(\xi)}I \right] u_0 \| \\
  &\le (1+M)K\frac{b_I}{a_I} \left(\frac{2}{a_I}\right)^\alpha \e^{-a_I t \cosh{\xi}-\alpha |\xi|} \|A^{\alpha}u_0\|,\; \xi \in \R, \; t\ge 0.
\end{split}
\]
The part responsible for the nonlocal condition in (\ref{par-int}),
can be estimated in the following way
\[
\left|\left(1 +\sum\limits_{k=1}^m \alpha_k
\e^{-z(\xi)t_k}\right)^{-1}\right| \le \left(1 -\sum\limits_{k=1}^m
|\alpha_k| \e^{-a_I\cosh(\xi)t_k}\right)^{-1}
\]
\cor{
\[
\le \left(1 -\sum\limits_{k=1}^m |\alpha_k|
\e^{-\rho_1t_k}\right)^{-1}\equiv Q^{-1},
\]
}
where $a_I\ge \rho_1$.

Thus, we obtain the following estimate for $\cF(t,\xi)$:
\begin{equation}\label{oc-pid}
    \|\cF(t,\xi) \|\le Q (1+M)K\frac{b_I}{a_I} \left(\frac{2}{a_I}\right)^\alpha \e^{-a_I t \cosh{\xi}-\alpha |\xi|}
    \|A^{\alpha}u_0\|,\; \xi \in \R, \; t\ge 0.
\end{equation}

The next step toward a numerical algorithm  is an approximation of (\ref{par-int}) by the efficient quadrature
formula. For this purpose we need to estimate the width of a strip around the real axis where the function $\cF(t,\xi)$
permit analytical extension (with respect to $\xi$). After changing $\xi$ to $\xi+i \nu$ the integration hyperbola
$\Gamma_I$ will be translated into the parametric hyperbola set
\begin{equation*}\label{par-hyp}
\begin{split}
\Gamma(\nu)&=\{z(w)=a_I \cosh{(\xi+i \nu)}-ib_I \sinh{(\xi+ i \nu)}: \; \xi\in(-\infty,\infty)\} \\
&=\{z(w)=a(\nu) \cosh{\xi}-ib(\nu) \sinh{\xi}: \; \xi\in(-\infty,\infty)\},
\end{split}
\end{equation*}
with
\begin{equation*}\label{na12}
\begin{split}
&a(\nu)= a_I \cos{\nu}+b_I
\sin{\nu}=\sqrt{a_I^2+b_I^2}\sin{(\nu+\phi/2)}, \\
& b(\nu)= b_I \cos{\nu}-a_I
\sin{\nu}=\sqrt{a_I^2+b_I^2}\cos{(\nu+\phi/2)},\\
& \cos{\frac{\phi}{2}}=\frac{b_I}{\sqrt{a_I^2+b_I^2}}, \;
\sin{\frac{\phi}{2}}=\frac{a_I}{\sqrt{a_I^2+b_I^2}} \; .
\end{split}
\end{equation*}

The analyticity of the function $\cF(t,\xi+i \nu),$ in the strip
\begin{equation*}\label{na16}
\begin{split}
&D_{d_1}=\{(\xi, \nu):\xi \in (-\infty, \infty), |\nu|<d_1/2\},
\end{split}
\end{equation*}
with some $d_1$ could be violated if the resolvent or the part related to the nonlocal condition become unbounded. To
avoid this we have to choose $d_1$ {in a way} that for $\nu \in (-d_1/2,d_1/2)$ the hyperbola set $\Gamma(\nu)$ remains
in the right half-plane of the complex plane. For $\nu=-d_1/2$ the corresponding hyperbola \cor{is} going {through the point
$(\rho_1,0)$, for some $0 \le \rho_1 <\rho_0$}. For $\nu=d_1/2$ it coincides with the spectral hyperbola and therefore
for all $\nu \in (-d_1/2,d_1/2)$ the set {$\Gamma(\nu)$} does not intersect the spectral sector. This fact justifies
the choice the hyperbola $\Gamma(0)=\Gamma_I$ as the integration path.

The requirements above imply the following system of equations
\begin{equation*}\label{na13}
  \begin{cases}
    a_I \cos{(d_1/2)}+b_I \sin{(d_1/2)}=\rho_0 ,   \\
   b_I \cos{(d_1/2)}-a_I \sin{(d_1/2)}=b_0=\rho_0\tan{\varphi} ,  \\
   a_I \cos{(-d_1/2)}+b_I \sin{(-d_1/2)}=\rho_1 ,
  \end{cases}
\end{equation*}
it lead us to the next system
\begin{equation*}\label{na14}
  \begin{cases}
   a_I=\rho_0 \cos{(d_1/2)} -b_0 \sin{(d_1/2)},  \\
   b_I=\rho_0 \sin{(d_1/2)} +b_0 \cos{(d_1/2)}, \\
   2 a_I \cos{(d_1/2)}=\rho_0+ \rho_1.
  \end{cases}
\end{equation*}
Eliminating $a_I$ from the first and the third equations we get
\[
\rho_0 \cos{d_1}- b_0 \sin{d_1}=\rho_1,
\]
\[
\cos(d_1+\varphi)=\frac{\rho_1}{\sqrt{\rho_0^2+b_0^2}},
\]
i.e.
\begin{equation}\label{shyr-sm}
d_1=\arccos{\left(\frac{\rho_1}{\sqrt{\rho_0^2+b_0^2}}\right)} -\varphi ,
\end{equation}
with $\cos{\varphi}=\frac{\rho_0}{\sqrt{\rho_0^2 +b_0^2}},$ $\sin{\varphi}=\frac{b_0}{\sqrt{\rho_0^2+b_0^2}}.$ Thus,
for $a_I$, $b_I$ we receive
\begin{equation}\label{na15}
\begin{split}
a_I&=\sqrt{\rho_0^2+b_0^2} \cos{\left(\frac{d_1}{2} +\varphi\right)} \\
& = \rho_0 \frac{\cos{\left(\frac{d_1}{2} +\varphi\right)}}{\cos{\varphi}} =\rho_0 \frac{\cos{\left(\arccos\left(\frac{\rho_1}{\sqrt{\rho_0^2
+b_0^2}}\right)/2  +\varphi/2 \right)}}{\cos{\varphi}}, \\
b_I&=\sqrt{\rho_0^2+b_0^2} \sin{\left(\frac{d_1}{2} +\varphi\right)} \\
&=\rho_0 \frac{\cos{\left(\frac{d_1}{2} +\varphi\right)}}{\cos{\varphi}} =\rho_0 \frac{\cos{\left(\arccos\left(\frac{\rho_1}{\sqrt{\rho_0^2
+b_0^2}}\right)/2  +\varphi/2 \right)}}{\cos{\varphi}}.
\end{split}
\end{equation}
{For $a_I$ and $b_I$ defined as above the vector valued function
$\cF(t,w)$ is analytic in the strip $D_{d_1}$ with respect to
$w=\xi+i \nu$  for any $t \ge 0$.} Note, that for $\rho_1=0$ we have
$d_1=\pi/2-\varphi$ as in \cite{gm5}.

Taking into account (\ref{na15}) we similarly can write equations for $a(\nu), b(\nu)$ on the whole interval
$-\frac{d_1}{2}\le\nu \le \frac{d_1}{2}$
\[
\begin{split}
a(\nu) &= a_I \cos{\nu}+ b_I\sin{\nu} =\sqrt{\rho_0^2+b_0^2}\cos{\left(\frac{d_1}{2} +\varphi\right)} \cos(\nu) \\
 &+\sqrt{\rho_0^2+b_0^2}\sin{\left(\frac{d_1}{2}+\varphi\right)} \sin(\nu) =\sqrt{\rho_0^2+b_0^2}\cos{\left(\frac{d_1}{2}+\varphi-\nu\right) } , \\
 b(\nu) &= b_I \cos{\nu}-a_I \sin{\nu} =\sqrt{\rho_0^2+b_0^2}\sin{\left(\frac{d_1}{2} +\varphi\right)} \cos(\nu) \\
 &-\sqrt{\rho_0^2+b_0^2}\cos{\left(\frac{d_1}{2}+\varphi\right)} \sin(\nu) =\sqrt{\rho_0^2+b_0^2}\sin{\left(\frac{d_1}{2}+\varphi -\nu\right)} ,
\end{split}
\]
\[
 \rho_1\le a(\nu)\le \rho_0, \quad b_0\le b(\nu)\le \sqrt{b_0^2+\rho_0^2 -\rho_1^2},
\]
with $d_1,$ defined by (\ref{shyr-sm}).

For the part responsible for the nonlocal condition we have
\[
\left|\left(1 +\sum\limits_{i=1}^n \alpha_i
\e^{-z(\xi,\nu)t_i}\right)^{-1}\right| \le \left(1
-\sum\limits_{i=1}^n |\alpha_i| |\e^{-z(\xi,\nu)t_i} |\right)^{-1}
\]
\[
\le \left(1 -\sum\limits_{i=1}^n |\alpha_i|
\e^{-a(\nu)\cosh(\xi)t_i}\right)^{-1}
\]
\[
\le \left(1 -\sum\limits_{i=1}^n |\alpha_i| \e^{-a(\nu)
t_i}\right)^{-1} < \left(1 -\sum\limits_{i=1}^n |\alpha_i|
\e^{-\rho_1t_i}\right)^{-1}<q^{-1}=Q,
\]
where $z(\xi,\nu)=a(\nu)\cosh(\xi)-ib(\nu)\sinh(\xi).$

Now, for $w \in D_{d_1}$ we get the estimate
\[
\|\cF(t,w)\| \le \e^{-a(\nu) t \cosh{\xi}}\frac{(1+M)QK\sqrt{a^2(\nu) \sinh^2{\xi}+b^2(\nu) \cosh^2{\xi}}}{(a^2(\nu) \cosh^2{\xi}
+b^2(\nu) \sinh^2{\xi})^{(1+\alpha)/2}}\|A^{\alpha}u_0\|
\]
\[
\le (1+M)QK\frac{b(\nu)}{a(\nu)}\frac{ \e^{-a(\nu) t \cosh{\xi}}}{(a^2(\nu) \cosh^2{\xi}+b^2(\nu) \sinh^2{\xi})^{(\alpha/2)}}\|A^{\alpha}u_0\|
\]
\[
\le (1+M)Q K\frac{b(\nu)}{a(\nu)} \left(\frac{2}{a(\nu)}\right)^\alpha \e^{-a(\nu) t \cosh{\xi}-\alpha |\xi|} \|A^{\alpha}u_0\|
\]
\[
\le (1+M)QK \tan{\left(\frac{d_1}{2} +\varphi-\nu \right)} \left(\frac{2\cos{\varphi}}{\rho_0 \cos{\left(\frac{d_1}{2}
+\varphi-\nu \right)}}\right)^\alpha  \e^{-\alpha |\xi|}\|A^{\alpha}u_0\|,
\]
\[
\forall w \in D_d.
\]

Similarly to \cite{stenger},  we introduce the space $ {\bf
H}^p(D_d), \; 1 \le p \le \infty$  of all vector-valued functions $\cal{F}$ analytic in
the strip
\[
D_d=\{z \in \mathbb{C}: - \infty < \Re z < \infty, |\Im z|<d \},
\]
equipped by the norm
\[
\| {\cal{F}} \|_{{\bf H}^p (D_d)}=
  \begin{cases}
    \lim_{\epsilon \to 0}(\int_{\partial D_d(\epsilon)}\|{\cal F}(z)\|^p |dz|)^{1/p} & \text{if $1 \le p < \infty$}, \\
   \lim_{\epsilon \to 0}\sup_{z \in \partial D_d(\epsilon)}\|{\cal F}(z)\| & \text{if
   $p=\infty$},
  \end{cases}
\]
where
\[
D_d(\epsilon)=\{z \in \mathbb{C}: | \text{Re}(z)| < 1/\epsilon,
|\text{Im}(z)|<d(1-\epsilon)\}
\]
and $\partial D_d(\epsilon)$ is the boundary of $ D_d(\epsilon)$.

Taking into account that the integrals over the vertical sides of
the rectangle $D_{d_1}(\epsilon)$ vanish as $\epsilon \to 0$, the
above estimate for $\|\cF(t,w)\|$ implies
\begin{equation}\label{H-norm}
\begin{split}
& \|\cF(t,\cdot)\|_{{\bf H}^1(D_{d_1})}\le  \|A^{\alpha}u_0\|[C_{-}(\varphi,\alpha) \\
 &+C_{+}(\varphi,\alpha)] \int_{-\infty}^\infty e^{-\alpha |\xi|} d \xi =C(\varphi,\alpha) \|A^{\alpha}u_0\|
 \end{split}
\end{equation}
with
\[
\begin{split}
&C(\varphi,\alpha)=\frac{2}{\alpha}[C_{+}(\varphi,\alpha) +C_{-}(\varphi,\alpha)], \\
&C_{\pm}(\varphi,\alpha)= (1+M)QK \tan{\left(\frac{d_1}{2} +\varphi \pm \frac{d_1}{2} \right)}
\left(\frac{2\cos{\varphi}}{\rho_0 \cos{\left(\frac{d_1}{2} +\varphi \pm \frac{d_1}{2} \right)}}\right)^\alpha .
\end{split}
\]
Note that the influence of both the smoothness parameter of $u_0$ given by $\alpha$ and of the spectral characteristics
of the operator $A$ \cor{given by $\varphi$ and $\rho_0$} is accounted by that fact, that the constant $C(\varphi,\alpha)$ from
(\ref{na15}) tends to $\infty$ if $\alpha \to 0,$ $\varphi \to \pi/2$ or $\rho_1 \to 0$ (in this case due to
(\ref{shyr-sm}) $d_1 \to \frac{\pi}{2}-\varphi $).

We approximate  integral (\ref{par-int}) by the following
Sinc-quadrature \cite{stenger, gm5}:
\begin{equation}\label{h-nab}
u_{h,N}(t)=\frac{h}{2 \pi i}\sum_{k=-N}^{N}\cF(t,z(kh)),
\end{equation}
{with an error }
\[
\|\eta_N(\cF,h)\|=\|u_h(t)-u_{h,N}(t)\|
\]
\[
\le \|u_h(t)-\frac{h}{2 \pi i}\sum_{k=-\infty}^{\infty}\cF(t,z(kh))\| +\|\frac{h}{2 \pi i}\sum_{|k|>N}\cF(t,z(kh))\|
\]
\[
\le \frac{1}{2 \pi}\frac{e^{-\pi d_1/h}}{2 \sinh{(\pi d_1/h)}}\|\cF\|_{{\bf H}^1(D_{d_1})}
\]
\[
+\frac{C(\varphi,\alpha)h \|A^{\alpha}u_0\| }{2\pi} \sum_{k=N+1}^{\infty}\text{exp}{[-a_I t \cosh{(kh)}-\alpha k h]}
\]
\[
\le \frac{c \|A^{\alpha}u_0\|}{\alpha}\left\{\frac{e^{-\pi d_1/h}}{\sinh{(\pi d_1/h)}} +\text{exp}{[-a_I t \cosh{((N+1)h)}-\alpha (N+1)h]}\right\},
\]
where the constant $c$ does not depend on $h,N,t.$ {Equalizing } the
both exponentials for $t=0$ implies
\[
 \frac{2 \pi d_1}{h}=\alpha(N+1)h ,
\]
or after the transformation
\begin{equation}\label{st-size}
 h=\sqrt{\frac{2 \pi d_1}{\alpha(N+1)}}.
\end{equation}
With this step-size the following error estimate holds true
\begin{equation}\label{na24}
\|\eta_N(\cF,h)\|\le \frac{c}{\alpha}\text{exp}{\left(-\sqrt{\frac{\pi d_1 \alpha}{2}(N+1)}\right)} \|A^{\alpha}u_0\|,
\end{equation}
with a constant $c$ independent of $t,N.$ In the case $t>0$ the first \cor{summand in} the argument of $\text{exp}{[-a_I t
\cosh{((N+1)h)}-\alpha (N+1) h]}$ from the estimate for $\|\eta_N(\cF,h)\|$ contributes mainly to the error order.
Setting in this case $h=c_1 \ln{N}/N$ with some positive constant $c_1$ we remain, asymptotically for a fixed $t,$ with
an error
\begin{equation}\label{na240}
\|\eta_N(\cF,h)\|\le c\left[e^{- \pi d_1 N/(c_1 \ln{N})} +e^{-c_1 a_I t N/2-c_1 \alpha \ln{N} }\right]\|A^{\alpha}u_0\|,
\end{equation}
where $c$ is a positive constant. Thus, we have proved the
following result.

\begin{theorem}\label{hom-conv}
Let $A$ be a densely defined strongly positive operator  and $u_0
\in D(A^{\alpha}),$ $\alpha \in (0,1),$ \cor{then the Sinc-quadrature}
(\ref{h-nab})  represents an approximate solution of the homogeneous
nonlocal value problem (\ref{pr-st}) (i.e. the case when $f(t)\equiv
0$) and \cor{possesses an exponential convergence rate which is uniform with respect to $t\ge 0$ and is of the order} ${\cal O}(e^{-c\sqrt{N}})$
uniformly in $t\ge 0$ provided that $h={\cal O}(1/\sqrt{N})$
(estimate (\ref{na24})) and of the order ${\cal
O}\left(\text{max}\left\{\e^{- \pi d N/(c_1 \ln{N})}, \right.
\right. $ $\left. \left.  \e^{-c_1 a_I t N/2-c_1 \alpha \ln{N} }
\right\}\right)$ for each fixed $t > 0$ provided that $h=c_1
\ln{N}/N$ (estimate (\ref{na240})).
\end{theorem}


\section{Inhomogeneous nonlocal problem}\label{s-inh-pr}

In this section we consider the particular solution
(\ref{bp1IntRed}) of inhomogeneous problem (\ref{pr-st}), i.e. with
$f(t)\neq 0$.

Let us rewrite formula (\ref{bp1IntRed}) in the form
\begin{equation}\label{ih1-nab-3}
\begin{split}
&u_{ih}(t)=u_{1,ih}(t)+u_{2,ih}(t),
\end{split}
\end{equation}
with
\begin{equation}\label{ih1-nab-2}
\begin{split}
&u_{1,ih}(t)=\int_0^t{e^{-A(t-\tau)}f(\tau)}d\tau, \quad
u_{2,ih}(t)=-\sum_{j=1}^m \alpha_ju_{2,ih,j}(t),
\end{split}
\end{equation}
where
\cor{
\begin{equation}\label{ih1-nab-1}
u_{2,ih,j}(t)= \int_0^{t_j} B^{-1} \e^{-A (t+t_j-\tau)} f(\tau) d\tau.
\end{equation}
}
We approximate the term {$u_{1,ih}(t)$} by the algorithm proposed in
\cite{gm5}:
\begin{equation}\label{ih1-nab}
\begin{split}
u_{1,ih}(t)\approx u_{1,N}(t)=&\frac{h}{2 \pi
i}\sum_{k=-N}^{N}z^{\prime}(kh)[(z(kh)I-A)^{-1}-\frac{1}{z(kh)}I]\\
& \times h\sum_{p=-N}^{N}\mu_{k,p}(t)f(\omega_p(t)),
\end{split}
\end{equation}
where
\[
\mu_{k,p}(t)=\frac{t}{2}\text{exp}\{-\frac{t}{2} z(kh)[1-\tanh{(ph)}]\}/\cosh^2{(ph)},
\]
\[
\omega_p(t)=\frac{t}{2}[1+\tanh{(ph)}], \; h={\cal
O}\left(1/\sqrt{N}\right),
\]
\[
z(\xi)=a_I \cosh{\xi}-i b_I \sinh{\xi},\; z'(\xi)=a_I \sinh{\xi}-i
b_I \cosh{\xi}.
\]

The next theorem (see \cite{gm5}) characterizes the error of this algorithm.
\begin{theorem}\label{inh1-conv}
Let $A$ be a densely defined strongly positive operator with
spectral characteristics $\rho_0,$ $\varphi$ and  a right
hand side $f(t) \in D(A^{\alpha}),$ $\alpha>0$ for $t \in
[0,\infty]$ can be analytically extended into the sector
$\Sigma_f=\{\rho \e^{i \theta_1 }: \; \rho \in [0,\infty], \;
|\theta_1|<\varphi\}$ where the estimate
\begin{equation}\label{analie24000}
  \|A^{\alpha}f(w)\|\le c_\alpha e^{-\delta_\alpha |\text{Re}\; w|}, \; w \in \Sigma_f
\end{equation}
with  $\delta_\alpha \in (0,\sqrt{2}\rho_0]$ holds, then algorithm
(\ref{ih1-nab}) converges to the solution of (\ref{pr-st}) with the error estimate
\begin{equation}\label{ie250}
\|\cE_N(t)\|=\|u_{1,ih}(t)-u_{1,N}(t)\|\le c \e^{-c_1\sqrt{N}}
\end{equation}
uniformly in $t$ with positive constants $c,$ $c_1$ depending on
$\alpha,$ $\varphi,$ $\rho_0$ and independent of $N.$
\end{theorem}

Let us construct an exponentially convergent approximation to the
term $u_{2,ih}$. {Using} the representation of the operator
functions by means of the Dunford-Cauchy integral for the j-th
summand of $u_{2,ih}$ we get
\begin{equation}\label{rep-ih2}
\begin{split}
 u_{2,ih,j}(t)=&\int_0^{t_j} \frac{1}{2 \pi i}\int_{\Gamma_I} \e^{-z(t+t_j-s)}B^{-1}(z)[(zI-A)^{-1}-\frac{1}{z}I]f(s)dz ds\\
=& \frac{1}{2 \pi i}\int_{\Gamma_I} \e^{-z(\xi)t} B^{-1}(z) \left[(z(\xi)I-A)^{-1}-\frac{1}{z(\xi)}I\right] \\
 &\times \int_0^{t_j} \e^{-z(\xi)(t_j-s)}f(s) ds z^{\prime}(\xi) d\xi, \\
z(\xi)=& a_I\cosh{\xi}- i b_I \sinh{\xi}.
\end{split}
\end{equation}
Replacing here the first integral by the Sinc-quadrature with
$h={\mathcal O}\left(N^{-1/2}\right)$ we obtain
\begin{equation}\label{ie4}
\begin{split}
u_{2,ih,j}(t)\approx u_{2,N}(t)= &\frac{h}{2 \pi i}\sum_{k=-N}^N \e^{-z(kh)t} z^{\prime}(k h) B^{-1}(z(kh))\\
& \times \left[(z(kh)I-A)^{-1}-\frac{1}{z(kh)}I\right] f_{k,j},
\end{split}
\end{equation}
where
\begin{equation}\label{ie5}
f_{k,j}=\int_0^{t_j} \e^{-z(kh)(t_j-s)}f(s) ds,\; k=-N,...,N,\; j=\overline{1,n}.
\end{equation}

In order to construct an exponentially convergent quadrature for these integrals we change the variables by
\begin{equation}\label{ie6}
s=\frac{t_j}{2}(1+\tanh{\xi}),
\end{equation}
and  obtain the improper integrals
\begin{equation}\label{ie7}
f_{k,j}=\int_{-\infty}^{\infty}{\cF}_k(t_j,\xi) d\xi,
\end{equation}
with
\begin{equation}\label{ie8}
{\cF}_k(t_j,\xi)=\frac{t_j}{2\cosh^2{\xi}}\text{exp}{[-z(kh)t_j(1-\tanh{\xi})/2]}f(t_j(1+\tanh{\xi})/2).
\end{equation}
Note that 
equation (\ref{ie6}) represents the conformal mapping {$w=\psi(z)=t_j[1+\tanh{z}]/2,$}
$z=\phi(w)=\frac{1}{2}\ln{\frac{t_j-w}{w}}$ of the strip $D_{\nu}$ onto the eye-shaped domain domain $\cA_{\nu}$
\cite{stenger}. On the real axes the integrand can be estimated by
\begin{equation}\label{ie9}
\begin{split}
\|{\cF}_k(t_j,\xi)\| \le &\frac{t_j}{2\cosh^2{\xi}} \text{exp}{[-a_I\cosh{(kh)}t_j(1 -\tanh{\xi})/2]} \\
& \qquad \quad \times \|f(t_j(1+\tanh{\xi})/2)\| \\
\le &2 t_j \e^{-2|\xi|}\|f(t_j(1+\tanh{\xi})/2)\|.
\end{split}
\end{equation}

The following result from \cite{gm5} is needed to justify  a
Sinc-quadrature for (\ref{ie7}):

\begin{lemma}\label{analytic-f}
Let the \cor{right-hand} side $f(t)$ in (\ref{pr-st}) for $t \in [0,\infty]$ can be analytically extended into
the sector $\Sigma_f=\{\rho e^{i \theta_1 }: \; \rho \in [0,\infty], \; |\theta_1|<\varphi\}$ and for all complex $w \in \Sigma_f$ we have
\begin{equation}\label{ie90}
  \|f(w)\|\le c e^{-\delta |\text{Re} \; w|},
\end{equation}
with  $\delta \in (0,\sqrt{2}\rho_0]$, then the integrand
${\cF}_k(t_j,\xi)$ can be analytically extended into a strip
$D_{d_1}, \; 0<d_1<\varphi/2$ and belongs to the class ${\bf
H}^1(D_{d_1})$ with respect to $\xi$, where $\rho_0,$ $\varphi$ \cor{are}
the spectral characteristics of $A.$
\end{lemma}

Under the assumptions of Lemma \ref{analytic-f} we can construct the
following exponentially convergent Sinc-approximation to $f_{k,j}:$
\begin{equation}\label{app-fkj}
f_{k,j} \approx f_{k,j,N} =h \sum_{p=-N}^{N} \mu_{k,p,j}f(\omega_{p,j}),
\end{equation}
where
\[
\mu_{k,p,j} =\frac{t_j}{2}\frac{\text{exp}\{-\frac{t_j}{2} z(kh) [1 -\tanh{(ph)}]\}}{\cosh^2{(ph)}},
\]
\[
\omega_{p,j} =\frac{t_j}{2}[1+\tanh{(ph)}], \quad h={\cal O}(1/\sqrt{N}),
\]
\[
z(\xi)=a_I \cosh{\xi}-i b_I \sinh{\xi}.
\]

Substituting  (\ref{app-fkj}) into (\ref{ie4}) we get the following
algorithm to compute an {approximation} $u_{2,j,N}$ to $u_{2,ih,j}$
\begin{equation}\label{ie4ff}
\begin{split}
u_{2,ih,j}(t)\approx u_{2,j,N}&(t)= \frac{h}{2 \pi i}\sum_{k=-N}^N \e^{-z(kh)t} z^{\prime}(k h) B^{-1}(z(kh))\\
& \times \left[(z(kh)I-A)^{-1}-\frac{1}{z(kh)}I\right] h \sum_{p=-N}^{N} \mu_{k,p,j}f(\omega_{p,j}).
\end{split}
\end{equation}

We represent the error in the form
\begin{equation}\label{ie26}
\cE_N (t) =u_{2,ih,j}(t) -u_{2,j,N}(t)=r_{1,N}(t)+r_{2,N}(t),
\end{equation}
where
\begin{equation}\label{ie27}
\begin{split}
&r_{1,N}(t)= u_{2,ih,j}(t) -u_{2,N}(t), \\
&r_{2,N}(t)=u_{2,N}(t) -u_{2,j,N}(t).
\end{split}
\end{equation}
Using  estimate (\ref{na24})(see also Theorem \ref{hom-conv})  we
obtain the estimate for $r_{1,N}(t)$
\begin{equation}\label{ie28}
\begin{split}
\|r_{1,N}(t)\| =&\left\|\int_0^{t_j} \left\{ \frac{1}{2 \pi i} \int_{-\infty}^\infty F_A(t+t_j-s,\xi)d \xi \right. \right. \\
 &\left. \left. -\frac{h}{2 \pi i}\sum_{k=-N}^N F_A(t+t_j-s,kh) \right\}f(s)ds \right\| \\
\le &\frac{c}{\alpha}\text{exp}{\left(-\sqrt{\frac{\pi d \alpha}{2}(N+1)} \right)} \int_0^{t_j} \|A^{\alpha}f(s)\|ds,
\end{split}
\end{equation}
where $F_A(t,\xi)$ is the operator defined in section
\ref{s-hom-pr}. Due to (\ref{fp11}) for $m=0$ we have for the error
$r_{2,N}$
\begin{equation}\label{ie29}
\begin{split}
\|r_{2,N}(t)\| =&\left\| \frac{h}{2 \pi i} \sum_{k=-N}^N \e^{-z(kh)t} z^{\prime}(kh) B(z(kh)) \right. \\
 &\left. \times\left[ (z(kh)I-A)^{-1}-\frac{1}{z(kh)}I \right] R_{k,j} \right \|\\
\le &\frac{h(1+M)QK}{2 \pi}\sum_{k=-N}^N \frac{|\e^{-z(kh)t}
z^{\prime}(kh)|}{|z(kh)|^{1+\alpha}}\|A^{\alpha}R_{k,j} \|,
\end{split}
\end{equation}
where
\[
R_{k,j}=f_{k,j}-f_{k,j,N}.
\]
The estimate (\ref{ie9}) yields
\begin{equation}\label{ie30}
 \|A^\alpha \cF_k(t_j,\xi)\|\le 2t_j \e^{-2|\xi|} \|A^\alpha f(\frac{t_j}{2}(1+\tanh{\xi}))\|.
\end{equation}
Due to Lemma \ref{analytic-f} the assumption $\|A^{\alpha}f(w)\|\le
c_\alpha e^{-\delta_\alpha |\text{Re} \; w|}$ $\forall w \in
\Sigma_f$ guarantees that $A^{\alpha}f(w) \in {\bf H}^1(D_{d_1})$
and $A^{\alpha}\cF_k(t_j,w) \in {\bf H}^1(D_{d_1}).$ These two condition turns us to the situation as in proposition of Theorem 3.2.1, p.144 from
\cite{stenger} with $A^{\alpha}f(w)$ instead of $f$ which implies
\[
\|A^{\alpha}R_{k,j} \| =\| A^{\alpha}(f_{k,j} -f_{k,j,N})\|
\]
\[
=\left\| \int_{-\infty}^\infty A^\alpha \cF_k (t_j,\xi)d \xi -h\sum_{k=-\infty}^{\infty}A^\alpha \cF_k(t_j,kh)\right\|
+ \left\| h\sum_{|k|>N}A^\alpha \cF_k (t_j,kh) \right\|
\]
\[
\le \frac{\e^{-\pi d_1/h}}{2 \sinh{(\pi d_1/h)}} \|\cF_k(t_j,w)
\|_{{\bf H}^1(D_{d_1})}
\]
\[
+h \sum_{|k|>N}2t_j \e^{-2|kh|}\|A^\alpha f(\frac{t_j}{2}(1+\tanh{kh}))\|
\]
\[
\le c \e^{-2\pi d_1/h}\|A^{\alpha}f(t_j,w)\|_{{\bf H}^1(D_{d_1})}
\]
\[
+h \sum_{|k|>N}2t_j \e^{-2|kh|}c_\alpha \text{exp} \left\{-\delta_\alpha \frac{t_j}{2}(1+\tanh{kh}) \right\} ,
\]
i.e., we have
\begin{equation}\label{ie31}
\|A^{\alpha}R_{k,j} \| \le c \e^{-c_1\sqrt{N}},
\end{equation}
where positive constants $c_{\alpha}, \delta_{\alpha},c,c_1$ do
not depend on $t,$ $N,$ $k$. Now, (\ref{ie29}) takes the form
\begin{equation}\label{ie32}
\begin{split}
\|r_{2,N}(t)\| &=\frac{h}{2 \pi i} \sum_{k=-N}^N \e^{-z(kh)t} z^{\prime}(kh)B(z(kh)) \\
&\times \left[ (z(kh)I-A)^{-1} -\frac{1}{z(kh)}I\right] R_{k,j}\\
&\le  c \e^{-c_1\sqrt{N}}S_N(t),
\end{split}
\end{equation}
with $S_N(t) =\sum_{k=-N}^N h \frac{|\e^{-z(kh)t} z^{\prime}(kh)|}{|z(kh)|^{1+\alpha}}.$ Using the estimate (4.8) from \cite{gm5} and
\begin{equation}\label{ie3300}
\begin{split}
&|z(kh)|=\sqrt{a_I^2 \cosh^2{(kh)} +b_I^2 \sinh^2{(kh)}}\\
&\ge a_I\cosh{(kh)}\ge a_I \e^{|kh|}/2,
\end{split}
\end{equation}
the last sum can be estimated by
\begin{equation}\label{ie33}
\begin{split}
&|S_N(t) |\le \frac{c}{\sqrt{N}}
\sum_{k=-N}^{N}\e^{-\alpha|k/\sqrt{N}|}\le c
\int_{-\sqrt{N}}^{\sqrt{N}} \e^{-\alpha t} dt\le c/\alpha  \;
\forall \; t\in [0,\infty).
\end{split}
\end{equation}
Taking into account (\ref{ie31}), (\ref{ie33}) we {deduce} from
(\ref{ie32})
\begin{equation}\label{ie34}
\|r_{2,N}(t)\|\le c \e^{-c_1\sqrt{N}}.
\end{equation}

The next assertion now follows  from (\ref{ie26}), (\ref{ie28}) and
(\ref{ie34}).

\begin{theorem}\label{main2}
Let the assumptions \cor{of theorem} \ref{inh1-conv} hold. Then
algorithm (\ref{ie4ff}) converges uniformly with respect to $t$ and moreover
the following error estimate holds true:
\begin{equation}\label{ie250-2}
\|\cE_N(t)\| =\|u_{2,ih,j}(t) -u_{2,j,N}(t)\|\le c \e^{-c_1\sqrt{N}},
\end{equation}
with positive constants $c,$ $c_1$ depend on $\alpha,$ $\varphi,$ $\rho_0$ and independent of $N.$
\end{theorem}

Now we can use approximation (\ref{ih1-nab}) for the every summand
and get:
\begin{equation}\label{inhnloc}
 u_{2,ih}(t)=\sum_{i=1}^n \alpha_i \int_0^{t_i} B^{-1}\e^{-A (t+t_i-\tau)} f(\tau) d\tau \approx \sum_{j=1}^n u_{2,j,N}(t)=u_{2,N}(t).
\end{equation}
Theorem \ref{main2} \cor{guarantees} that the error
this approach will be bounded by:
\begin{equation}\label{inh2-conv}
    \|u_{2,ih}(t)-u_{2,N}(t)\| \le \sum_{j=1}^m c \e^{-c_1\sqrt{N}} \le c_2 \e^{-c_1\sqrt{N}}.
\end{equation}

Thus, the approximations (\ref{h-nab}) together with (\ref{ih1-nab})
and (\ref{inhnloc}) represent an exponentially convergent algorithm
for the problem (\ref{pr-st}).


\section{Numerical examples}\label{num-ex}
\begin{example}
We consider the homogeneous problem (\ref{pr-st}) with the operator $A$ defined  by
\begin{equation}\label{oper-def}
\begin{split}
&D(A)=\{u(x) \in H^2(0,1): \; u(0)=u(1)=0\}, \; \\
&Au=-u^{\prime \prime}(x) \; \forall u \in D(A).
\end{split}
\end{equation}
The initial nonlocal condition reads as follows:
\[
 u(x,0)+0.5u(x,0.2)+0.3u(x,0.4)=(1+0.5\e^{-\pi^2 0.2}+0.3\e^{-\pi^2 0.4})\sin(\pi x),
\]
with $u_0=(1+0.5\e^{-\pi^2 0.2}+0.3\e^{-\pi^2 0.4})\sin(\pi x) \in
D(A)$. The exact solution of the problem is
$u(x,t)=\e^{-\pi^2t}\sin(\pi x).$ It is easy to find that
\[
(zI-A)^{-1}u_0=\left( z+ \frac{d^2}{dx^2} \right)^{-1}\sin(\pi x)=\frac{\sin(\pi x)}{z-\pi^2}.
\]

Calculations with the algorithm {(\ref{ie4ff})} above has been performed in Maple with $h=N^{-1/2}.$ The error at
$x=0.5, \; t=0.3$ is presented in Table \ref{tab1} and clearly exhibits an exponential decay \cor{according to the theoretical estimates}.
\begin{table}[ptbh]
\begin{center}
  \begin{tabular}{|c|c|}
    \hline
   N & $\varepsilon_N$ \\
   \hline
    4 &   .29857983847712589e-1\\
    8 &   .41823888073604986e-2 \\
    16 &  .11258594468208641e-2 \\
    32 &  .10042178166563831e-3 \\
    64 &  .28007158539828452e-5 \\
    128 & .2098826601399176e-7 \\
    256 & .1858929920173152e-10 \\
    512 & .856837124351510e-15 \\
        \hline
  \end{tabular}
\end{center}
\caption{The error for $x=0.5, \; t=0.3.$}\label{tab1}
\end{table}

Due to Theorem \ref{hom-conv} the error {should not be greater then}
$ \varepsilon_N={\cal O}\left(\e^{-c\sqrt{N}}\right).$ The constant
$c$ in the exponent can be estimated using the following
a-posteriori relation:
\[
c=\ln\left(\frac{\varepsilon_N}{\varepsilon_{2N}}\right)(\sqrt{2}-1)^{-1}N^{-1/2}= \ln\left(\mu_N \right)(\sqrt{2}-1)^{-1}N^{-1/2}.
\]
{The numerical results are presented in the Table \ref{tab-c} for
this estimation show that the constant can be estimated as $c\approx
1.5$  when $N\rightarrow \infty$. }
\begin{table}[ptbh]
\begin{center}
  \begin{tabular}{|c|c|}
    \hline
    N &          $c$ \\
  \hline
    4 &   2.372652515388745588587496\\
    8 &   1.120148732795449515627946 \\
    16 &  1.458741976765153165445005 \\
    32 &  1.527648924601130131250452 \\
    64 &  1.476794596387591759032900 \\
    128 & 1.499935011373075736075927 \\
    256 & 1.506597339081609844717370 \\
  \hline
  \end{tabular}
\end{center}
\caption{The estimate of $c$}\label{tab-c}
\end{table}

\end{example}

{The next example deals again with a homogeneous problem but in this
more realistic case the resolvent of $A$ on the element  $u_0$ \cor{cannot}
be calculated analytically. }
\begin{example}
We consider the homogeneous problem (\ref{pr-st}) with the operator
$A$ defined as in (\ref{oper-def}) and  with the following initial
nonlocal condition:
\[
 u(x,0)+ u(x,0.5)= x \ln(x),
\]
where $u_0=x\ln(x) \in A^{\alpha},$ $\alpha < 1/2.$ In this case the
resolvent can be represented using the Green function
\[
(zI-A)^{-1}u_0=\left( z+ \frac{d^2}{dx^2} \right)^{-1}x\ln(x)=\int_0^1 G(x,s) s\ln( s)ds,
\]
\[
G(x,s) =-\frac{1}{\sqrt{z} \sin(\sqrt{z})}
\begin{cases}
 \sin(x\sqrt{z}) \sin((1-s)\sqrt{z})   & x \le s, \\
 \sin(s\sqrt{z}) \sin((1-x)\sqrt{z}) & x \ge s
 \end{cases},
\]
where the integrals were computed by exponentially convergent
Sinc-quadrature (see e.g. \cite{stenger}) using  Maple. The results
for $x=0.5, \; t=0.3$ are presented in Table \ref{tab2}.
\begin{table}[ptbh]
\begin{center}
  \begin{tabular}{|c|c|}
    \hline
   N & $u(x,t)$ \\
   \hline
    4 &    -.241535790017043e-1\\
    8 &    -.228401191029108e-1\\
    16 &   -.194273285627507e-1 \\
    32 &   -.192905848633180e-1 \\
    64 &   -.192911920318628e-1 \\
    128 &  -.192907849909929e-1 \\
    256 &  -.192907820740651e-1 \\
        \hline
  \end{tabular}
\end{center}
\caption{Values of the solution $u(x,t)$ for $x=0.5, \; t=0.3$.}\label{tab2}
\end{table}

It can be easily seen that the number of stabilized digits increases
according to the theoretical prediction by Theorem \ref{hom-conv}.
\end{example}

\begin{example}
Let us  consider the inhomogeneous problem (\ref{pr-st}) with the same $A$ defined  by (\ref{oper-def}), and the
nonlocal condition
\[
 u(x,0)+0.5u(x,0.2)=(1+0.5\e^{0.2})\sin(\pi x),
\]
For $f(t,x)$ at the \cor{right-hand} side of the equation (\ref{pr-st}) we set
\[
f(x,t)=(1+\pi^2)\e^t\sin(\pi x).
\]
The exact solution of the problem is $u(x,t)=\e^{t}\sin(\pi x).$ We have used the algorithm defined by (\ref{h-nab}),
(\ref{ih1-nab}), (\ref{inhnloc}) and implemented in Maple.  For $x=0.5, \; t=0.3$  the results presented in Table
\ref{tab3} and are again in good agreement with the theoretical predictions.
\begin{table}[ptbh]
\begin{center}
  \begin{tabular}{|c|c|}
    \hline
   N & $\varepsilon_N$ \\
   \hline
    4 &   .202211483120243\\
    8 &   .726677678737409e-1 \\
    16 &  .138993889900620e-1 \\
    32 &  .143037059411419e-2 \\
    64 &  .554542099757830e-4 \\
    128 & .532640823981411e-6 \\
    256 & .730569324317506e-9 \\
    512 & .648376079810788e-13 \\
        \hline
  \end{tabular}
\end{center}
\caption{The error for $x=0.5, \; t=0.3.$}\label{tab3}
\end{table}
\end{example}


{\bf Acknowledgment}. The authors would like to acknowledge the
support provided by the Deutsche Forschungsgemeinschaft (DFG).

\end{document}